\documentclass[11 pt]{amsart}
\usepackage{amssymb}
\usepackage{amscd}
\usepackage{epsfig}
\newtheorem{theorem}{Theorem}
\newtheorem{corollary}[theorem]{Corollary}
\newtheorem{lemma}[theorem]{Lemma}
\newtheorem{proposition}[theorem]{Proposition}

\theoremstyle{definition}
\newtheorem{definition}{Definition}

\numberwithin{equation}{section}

\newcommand{\bc}{\mathbb{C}}

\newcommand{\bz}{\mathbb{Z}}
\newcommand{\br}{\mathbb{R}}

\newcommand{\cg}{\mathcal{G}}

\newcommand{\hk}{\hookrightarrow}
\newcommand{\bg}{\bigskip}
\newcommand{\med}{\medskip}
\newcommand{\la}{\longrightarrow}
\newcommand{\bfl}{\begin{flushleft}}
\newcommand{\efl}{\end{flushleft}}

\newcommand{\gln}{U(n)}
\newcommand{\cge}{\cg (e^*(\zeta))}

\newcommand{\lu}{LU(n)}
\newcommand{\lcn}{L\bc^n}

\newcommand{\lc}{C^\infty (S^1 , \bc)}
\newcommand{\lplus}{L_+\bc^n}
\newcommand{\glk}{U(n)}
\newcommand{\lgk}{L\glk}

\newcommand{\om}{\Omega}

\newcommand{\restrict}{\!\!\mid}

\title{Fourier Decompositions of Loop Bundles}
\author[R.L. Cohen]{Ralph L. Cohen}
 \address{Dept. of Mathematics \\
Stanford University\\
Stanford, California 94305}
 \email[Cohen]
{ralph@math.stanford.edu}
\thanks{The first author was partially supported by a grant from the NSF
} 
\author[A. Stacey]{Andrew Stacey}
\address{Department of Mathematics\\
Stanford University \\
Stanford, California 94305}
\email[Stacey]{astacey@math.stanford.edu}
 \date{\today}
  
\begin{document}

\begin{abstract}
    In this paper we investigate bundles whose structure group is the
loop group $\lu$.  These bundles are classified by maps to the loop
space of the classifying space, $LBU(n)$.  Our main result is to give
a necessary and sufficient criterion for there to exist a Fourier type
decomposition of such a bundle $\xi$.  This is essentially a
decomposition of $\xi$ as $\zeta \otimes L\bc$, where $\zeta$ is a
finite dimensional subbundle of $\xi$ and $L\bc$ is the space of
functions, $ \lc$.  The criterion is a reduction of the structure
group to the finite rank unitary group $U(n)$ viewed as the subgroup
of $\lu$ consisting of constant loops.  Next we study the case where
one starts with an $n$ dimensional bundle $\zeta \to M$ classified by
a map $f : M \to BU(n)$ from which one constructs a loop bundle
$L\zeta \to LM$ classified by $Lf : LM \to LBU(n)$.  The tangent
bundle of $LM$ is such a bundle. We then show how to twist such a
bundle by elements of the automorphism group of the pull back of
$\zeta$ over $LM$ via the map $LM \to M$ that evaluates a loop at a
basepoint. Given a connection on $\zeta$, we view the associated
parallel transport operator as an element of this gauge group and show
that twisting the loop bundle by such an operator satisfies the
criterion and admits a Fourier decomposition.
\end{abstract}
\maketitle

\section*{Introduction  }
 In this paper we study infinite dimensional vector bundles whose
structure group is the loop group, $\lu$.  We call such a bundle a
``rank $n$ loop bundle". Here $U(n)$ is the unitary group of rank $n$
and for any finite dimensional smooth manifold $M$, $LM$ denotes the
space of smooth loops, $LM = C^\infty (S^1, M)$.  A rank $n$ loop
bundle will have fibers isomorphic to the loop space $\lcn$.  Now
$\bc^n$-valued functions on the circle have Fourier expansions and one
of our goals in this paper is to give necessary and sufficient
conditions for there to exist such an expansion fiberwise on a loop
bundle.

The Fourier expansion of elements in $\lcn$ can be viewed as a map
$$
\Phi: \lcn \la \bc [\![z, z^{-1}]\!] \otimes \bc^n
$$
where $\bc [\![z, z^{-1}]\!]$ denotes the ring of formal power series
in a variable $z$ and its inverse.  Let $\bc [\![z ]\!] \subset \bc
[\![z, z^{-1}]\!]$ be the subalgebra consisting of power series that
only involve non-negative powers of $z$.  Recall that the inverse
image, $\Phi^{-1} (\bc [\![z ]\!] \otimes \bc^n)$, which we call
$L_+\bc^n$, consists of boundary values of holomorphic
maps $f : D^2 \to \bc^n$.  Notice that $\lplus \subset \lcn$ has the
following properties:

 \begin{enumerate} \label{fourier}
\item Consider the inner product on $L\bc^n$ defined by
 \begin{equation} \label{eq:ip}
\left( \alpha,  \beta\right) = \int_{S^1} \langle \alpha(t), \beta(t) \rangle dt.    
\end{equation}
 Then 
$\lplus \subset \lcn$   is an infinite dimensional, infinite codimensional subspace that is closed with respect to the topology induced by
the inner product.  Moreover,    there is an orthogonal decomposition
$$
L\bc^n = \lplus \oplus L_- \bc^n
$$
 where the space $L_- \bc^n$ consists of loops whose Fourier series only involve negative powers of $z$.
 \item $\lplus$ 
  is invariant under multiplication by $z$: 
$$z\lplus
\subset
\lplus.$$  Here we are  identifying the Laurent polynomial ring
$\bc [z, z^{-1}]
$ as a subalgebra of $\lc$, and hence $\lcn$ is a module over $\bc [z, z^{-1}]$.  Furthermore, the codimension
of $z\lplus$
 in
$\lplus$ is $n$. 
 \end{enumerate}

\med
Notice that there is a filtration of subspaces,
$$
\cdots \subset  \cdots \subset z^{-k}\lplus  \subset z^{-(k+1)}\lplus \subset   \cdots  L\bc^n 
$$ where the union $\bigcup_k z^{-k}\lplus$ is a dense subspace of $L\bc^n$.
We think of this filtration as   the ``Fourier decomposition" of $\lcn$  and   will define the existence of a  ``Fourier
decomposition" of a loop bundle 
$\xi$ as the existence of a  subbundle, $\xi_+ \subset \xi$,  which satisfies the above  two  properties in a fiberwise manner. 

\med
One of the main goals of this paper is to prove the following theorem.    It can  be
viewed as saying that Fourier decompositions of loop bundles are most often  not possible.   

\med
\begin{theorem}\label{main}  Let $\xi \to X$ be a rank $n$ loop bundle.  Then $\xi$ has a Fourier decomposition if
and only if the structure group of $\xi$ can be reduced to $U(n)$
viewed as the subgroup of the constant loops in $\lu$. 
\end{theorem}

\med
  Notice that since the classifying space
of the infinite dimensional group $\lu$ is the loop space of the classifying space, of $U(n)$, $LBU(n)$, a rank $n$ loop bundle over
a space
$X$ is classified by  a map $f : X \to LBU(n)$.  We call the ``underlying $n$-dimensional bundle'',  $U(\xi)$, the bundle classified
by the composition
$$
\begin{CD} X @>f >>  LBU(n) @>e  >> BU(n) \end{CD}
$$
where $e   : LY \to Y$ is the map that evaluates a loop at $1 \in S^1$.
Let $D^2$ be the two dimensional disk.  Restriction to the boundary defines a fibration between the mapping spaces
$\rho : Map(D^2, BU(n)) \to LBU(n)$.  

\med
The following is  essentially a restatement  of the above theorem, but it is quite useful in practice. 

\med
\begin{corollary}  Let $f : X \to LBU(n)$ classify a loop bundle $\xi \to X$.  Then $\xi$ has a  Fourier decomposition
if and only if there is a lift of $f$ to $Map(D^2, BU(n))$.
 \end{corollary}

The following is an immediate corollary of theorem 1  and  gives a description of any loop bundle that has a Fourier decomposition.  

\begin{corollary}\label{description}  Let $\xi \to X$ be a rank $n$ loop bundle  with underlying $n$-dimensional bundle 
$U(\xi) \to X$.   Then $\xi $ has a Fourier decomposition if and only if there is an isomorphism of loop bundles
$$
\xi   \cong L\bc \otimes U(\xi).
$$
 \end{corollary}

\med
Let $M$ be a simply connected manifold and $\zeta \to M$ an $n$-dimensional complex bundle classified by a map
$f_\zeta : M \to B\gln$.    Consider the induced rank
$n$   loop bundle over the loop space, $L\zeta \to LM$.    
 
   Let $\gamma \in LM$.  The fiber over $\gamma$  of $L\zeta$
is the space of sections of the pull back of $\zeta$ over the circle, 
$$
L\zeta_\gamma =  \Gamma_{S^1} (\gamma^*(\zeta)).
$$

Perhaps the most important example of these bundles is when $\zeta = TM$ is the tangent bundle of $M$.
The bundle $LTM$ is the infinite dimensional tangent bundle of $LM$.  The  tangent space over $\gamma$ is   the
space of  vector fields living over $\gamma$.

Assume that $M$ is simply connected, so that $LM$ is connected.  In this case, for any $\gamma \in LM$, the pull
back, $\gamma^* (\zeta)$  is trivial.  Choosing a trivialization gives us isomorphisms 
\begin{equation}
L\zeta_\gamma = \Gamma_{S^1} (\gamma^*(\zeta)) \cong L\bc \otimes \zeta_{\gamma (0)}.
\end{equation}
Therefore,  for each $\gamma \in LM$,  we have a Fourier decomposition of the fiber, $L\zeta_{|_\gamma}$, induced by
the subspace corresponding to $L_+\bc \otimes \zeta_{\gamma (0)}$ with respect to such a trivialization.  However, this will not
necessarily vary continuously over $LM$.

The
following corollary to the above theorem says that
$L\zeta$ rarely will have a global Fourier decomposition.

\begin{corollary}\label{omM}  If $L\zeta \to LM$ has a Fourier decomposition, then the map of based loop spaces,
$$\om f_\zeta : \om M \to \om B\gln \simeq \gln$$
is null homotopic. 
\end{corollary}

 \bf Remark.  \rm  
 The
homotopy type of the  map of based loop spaces,
$\om f_\zeta :
\om M
\to
\gln  $  can be obtained by taking the holonomy map of a connection on $\zeta$.  Therefore this corollary can be
interpreted as saying that in order for  $L\zeta$  to have a Fourier decomposition, $\zeta$ must admit a ``homotopy flat"
connection.  

The parallel transport operator induced by a connection on $\zeta \to M$ can be interpreted as an automorphism
of the pull - back bundle $e^*(\zeta) \to LM$, where, as above, $e : LM \to M$ is the evaluation map. Let $\cg
(e^*(\zeta)) $ be the gauge group of bundle automorphisms of $e^*(\zeta)$.    Our final result  shows how to deform
the loop bundle $L\zeta$ by such an automorphism  and examines its decomposability properties.   

\med
\begin{theorem} \label{gauge} Let $\zeta \to M$ be an $n$ dimensional bundle over a smooth, simply connected
manifold.  Then there is a natural rank $n$ loop bundle
$$L^\cg \zeta \to \cge \times LM$$ satisfying the following properties.
For  $\tau \in \cge$, let $L^\tau M$ denote the restriction of $L^\cg
\zeta$ to $\{\tau\} \times LM$; then,
\begin{enumerate}
\item  For  the identity element, $id \in \cge$,
$L^{id}\zeta = L\zeta \to LM.$ 
\item  For $\tau_\alpha$ the parallel transport operator of a connection $\alpha$  on $\zeta$,  the bundle $L^{\tau_\alpha} \zeta \to LM$
admits a natural isomorphism of loop bundles,
 $$
L^\tau \zeta \cong C^\infty (S^1, \bc)\otimes e^*\zeta 
$$
and hence it admits   Fourier decomposition.
\end{enumerate}
\end{theorem}

\bf Remark.  \rm Given a connection on $\zeta$ with parallel transport operator $\tau$, we refer to the resulting
loop bundle $L^\tau \zeta$ has the ``holonomy loop bundle" induced by this connection.  Notice that the
isomorphism type of these bundles is independent of the connection because the space of connections is connected.   
Hence all  parallel transport operators coming from connections lie in the same path component of $\cge$.

\med
This paper is organized as follows.  In section 1 we will  give careful definitions of loop bundles and their Fourier decompositions and prove
theorem 1 and its corollaries.  Theorem 5 will be proved in section 2. 

\med
This paper was motivated by studying the beautiful ideas contained in J. Morava's work \cite{morava}.  
The authors are grateful to   Morava for many hours of discussion about loop bundles, holonomy, and related topics.

\bg 
\section{  Fourier decompositions and a proof of theorem 1}

In this section we define a rank \(n\) loop bundle and the concept of
a Fourier decomposition on such a bundle.  We conclude by proving
theorem~\ref{main}.

\med
In what follows all of our spaces will be of the homotopy type of a $CW$ - complex of finite type.

\med

\begin{definition}
  a.   A \emph{rank \(n\) loop bundle} \(\xi \to X\)
is a vector bundle over \(X\) with fiber isomorphic to  \(L\bc^n\) and structure group
\(\lgk\).   The classifying space of such a bundle is the loop space $LBU(n)$. 

b. Let \(\xi \to X\) be a rank \(n\) loop bundle classified by a map $f_\xi : X \to LBU(n)$. Let $e : LBU(n) \to BU(n)$ be the map that
evaluates a loop at $1 \in S^1$.  The \emph{underlying \(n\)-dimensional vector bundle},
\(U(\xi) \to X\), is the bundle classified by the composition $e \circ f_\xi : X \to LBU(n) \to BU(n)$.  
\end{definition}

In the language of principal bundles, a rank \(n\) loop bundle
classified by a map $f_\xi : X \to LBU(n) $ is:
\[
\xi = f_\xi^*(LEU(n)) \times_{\lgk} L\bc^n
\]
where $EU(n) \to BU(n)$ is the  universal $U(n)$-bundle.
The  underlying vector bundle is:
\[
U(\xi) = f_\xi^*(LEU(n)) \times_{\lgk} \bc^n = f_\xi^*e^*(E\glk)
\times_{\glk} \bc^n
\]
where \(\lgk\) acts on $\bc^n$  via the map \(\gamma \cdot v =
\gamma(1)v\).

As mentioned in the introduction, an important class of   loop bundles arises in the following way: let
\(\zeta \to M\) be a finite dimensional vector bundle over a finite
dimensional manifold \(M\).  Let \(f_\zeta : M \to B\glk\) be a
classifying map for \(\zeta\).  The loop space of \(\zeta\)  has the
structure of a loop bundle over \(LM\): \(L \zeta \to LM\).  It is classified by the map
\(Lf_\zeta : LM \to LB\glk \simeq B\lgk\).  If \(\zeta
\to M\) is the tangent bundle of \(M\), \(L\zeta\) is the tangent
bundle of \(LM\).

Another interesting class of examples of loop bundles arises by taking a finite
dimensional vector bundle \(\zeta \to X\) and forming the fiberwise
tensor product with \(L\bc\).  In terms of classifying maps, this
corresponds to the map \(\left[X, B\glk\right] \to \left[X,
LB\glk\right]\) induced by the inclusion \(B\glk \to LB\glk\) as the
space of constant maps.

\med
Since the structure group of our loop bundles is $LU(n)$, any such bundle
has an inner product that  on each fiber is isomorphic to the inner product on $L\bc^n$ given by (\ref{eq:ip}). 
This  allows us to
consider the existence of orthogonal complements of subbundles.  Such
complements may not necessarily exist, even when the subbundle is
fiberwise closed for the pre-Hilbertian topology defined by the inner
product.  But when a complement does exist locally, it can be extended
globally using the unitary structure.

Using the properties  of the subspace $L_+\bc^n \subset L\bc^n$ described in
the introduction  as our fiberwise model, we now define the notion
of a ``\sl Fourier decomposition"  \rm of a loop bundle. 

 \med 
\begin{definition}
A \emph{Fourier decomposition} of a   rank \(n\) loop bundle \(\xi
\to X\) is a  subbundle $\psi \subset \xi$ satisfying the following properties:
\begin{enumerate}
\item $\psi$ has an orthogonal complement $\psi^\perp \subset \xi$ with
$ 
 \xi = \psi \oplus \psi^\perp .
$ 
\item $\psi$ is invariant under multiplication by $z \in \bc[z, z^{-1}]$,
  $  z\psi \subseteq \psi. $
Furthermore $z\psi $   has   codimension \(n\) in $\psi$. 
\end{enumerate}
\end{definition}

\med

\bf Remarks. \rm 
1. The notation for the decomposition implies that \(\psi\) is fibrewise
closed in the pre-Hilbertian topology.  

2. We will see in our proof of theorem \ref{main}  that with the existence of such a subbundle,  there exists a filtration
$$
\cdots \subset z^{-k}\psi \subset z^{-(k+1)}\psi \subset \cdots \subset \xi
$$
whose union $\bigcup_k z^{-k} \psi$ is a  fiberwise dense subbundle of $\xi$. 
This is the bundle theoretic analogue of the Fourier decomposition of $L\bc^n$. 

3.  A Fourier decomposition yields a \sl polarization \rm  of $\xi$ in the sense
of \cite{segpres}.  However a Fourier decomposition is much stronger than a polarization 
since the splitting of a polarized bundle at a particular fiber need  only be well defined
up to a finite dimensional ambiguity.  See \cite{segpres} for details. 

\med

We are now ready to give a proof of theorem 1, as stated in the introduction.  We begin with a lemma.
 
\begin{lemma}
\label{prop:orth}
Let \(W \oplus W^\perp\) be an orthogonal decomposition of \(L\bc^n\)
with the property that \(zW \subseteq W\) has codimension \(n\).  Then
the space \(W \cap zW^\perp\) is an \(n\) dimensional subspace of
\(L\bc^n\).  Furthermore, given a unitary isomorphism, \(\phi : \bc^n \to W \cap zW^\perp\), the following composition defines an
element of $LU(n)$:     
 $$
\begin{CD} 
 L\bc \otimes \bc^n @>1 \otimes \phi > \cong > L\bc \otimes (W \cap zW^\perp) @>\rho >>
L\bc^n
\end{CD}
$$
 Here $\rho$ is the unique map of $L\bc$-modules  extending the inclusion $W\cap zW^{\perp} \hk L\bc^n$.   
\end{lemma}

\begin{proof}
 
Consider the space \(W \cap zW^\perp\).  This lies inside \(W\) and is
orthogonal to \(zW\) in \(W\).  As the map \(L\bc^n \to L\bc^n\)
defined by \(\alpha \to z\alpha\) is a unitary isomorphism, \(L\bc^n\)
can be orthogonally decomposed as \(zW \oplus zW^\perp\).  From this
and the fact that \(zW \subseteq W\) has codimension $n$,  we deduce that there is an
orthogonal decomposition \(W = (W \cap zW^\perp) \oplus zW\).  Hence \(W
\cap zW^\perp\) is \(n\) dimensional.

 We continue with an argument adapted from \cite[8.3.2]{segpres}.  Let \(\{w_1,
\dotsc, w_n\}\) be an orthogonal basis for \(W \cap zW^\perp\).  Let
\((w_{jl})\) denote the sequence of Fourier coefficients for \(w_j\).
As each \(w_j\) is a smooth map \(S^1 \to \bc^n\), the sum \(\sum_{l
\in \bz} w_{jl} e^{il\theta}\) is absolutely convergent and sums to
\(w_j(\theta)\).  We consider the inner product in \(\bc^n\) of
\(w_j(\theta)\) and \(w_k(\theta)\):
\begin{align*}
\langle w_j(\theta), w_k(\theta)\rangle &= \left<  \sum_{l \in \bz} w_{jl}
e^{il\theta}, \sum_{m \in \bz} w_{km} e^{im\theta} \right>\\
&= \sum_{l,m \in \bz} \left< w_{jl}, w_{km} \right> e^{i(m-l)\theta} \\
&= \sum_{p \in \bz} \sum_{r \in \bz} \left< w_{jr}, w_{k(p+r)} \right>
e^{ip\theta} \\
&= \sum_{p \in \bz} \left(w_j, z^pw_k\right) e^{ip\theta} \\
&= \delta_{jk}
\end{align*}
the last identity uses the fact that the \(w_j\) are an orthonormal
basis for \(W \cap zW^\perp\) and that this space is orthogonal to
\(zW\) within which the \(zw_j\) lie.

Hence \(\{w_1(\theta), \dotsc, w_n(\theta)\}\) is an orthonormal basis
for \(\bc^n\) for each \(\theta \in S^1\).  Thus there is an element
\(\gamma \in \lu\) such that \(\gamma(\theta)e_j = w_j(\theta)\).  The
fact that the \(w_j\) are smooth demonstrates that \(\gamma\) is
smooth.  The map defined in the statement of the proposition is easily
seen to agree with \(\gamma\).
\end{proof}

 The proof of theorem~\ref{main}  can now be completed quickly.

\begin{proof}[Proof of theorem~\ref{main}]

Let $\xi$ have a Fourier decomposition $\psi \oplus \psi^\perp \cong \xi$.  Consider the subspace, 
  \(\psi \cap z\psi^\perp \subset   \xi\).  We first show that
  \(\psi \cap z\psi^\perp \to X\)  is a subbundle of $\xi$.  That is,   we   show
that it is a locally trivial fibration with fiber a vector space of
constant rank.  Consider \(V \subseteq X\) open over which the
decomposition \(\xi = \psi \oplus \psi^\perp\) can be trivialised;
i.e.~there is an isomorphism   \(\xi \restrict_V \to V \times L\bc^n\)
which carries \(\psi \restrict_V\) to a subspace \(W\) and
\(\psi^\perp \restrict_V\) to \(W^\perp\).  This map is   \(L\bc\)
equivariant and so \(zW \subseteq W\) with codimension \(n\).  We note
that \(\psi \cap z\psi^\perp\) corresponds to \(W \cap zW^\perp\)
which is an \(n\) dimensional vector space.  Thus \(\psi \cap
z\psi^\perp\) is a locally trivial fibration with fiber isomorphic to  \(\bc^n\).

By lemma 6,   the evaluation map, \(e : L\bc^n \to \bc^n\), carries \(W \cap
zW^\perp\) isomorphically onto \(\bc^n\). This implies that the composition
$$
\begin{CD}
\psi \cap z\psi^\perp \hk \xi @> e >> U(\xi)
\end{CD}
$$
is an isomorphism.   

Now notice  that the map of $L\bc$-module bundles  \( \rho : L\bc \otimes ( \psi \cap z\psi^\perp )   \to \xi\)
extending the inclusion $\psi \cap z\psi^\perp  \hk \xi$ corresponds on the fibers over \(V\) to:
$$
\begin{CD}
L\bc^n \to L\bc \otimes \bc^n \to L\bc \otimes (W \cap zW^\perp) @>\rho >>
L\bc^n
\end{CD}
$$
By lemma \ref{prop:orth} this is an isomorphism. This implies 
$$
\rho : L\bc \otimes ( \psi \cap z\psi^\perp )   \to \xi
$$
is an isomorphism of loop bundles.  Coupled with the isomorphism $e :  \psi \cap z\psi^\perp \cong U(\xi)$, we have an isomorphism of
loop bundles
$$
 L\bc \otimes U(\xi) \cong \xi.
$$
But, as observed earlier,   the left hand bundle is classified by the composition
$$
\begin{CD}
X @>f_\xi >> LBU(n) @>e>> BU(n) @>\iota >> LBU(n)
\end{CD}
$$
where $\iota : BU(n) \hk LBU(n)$ is the inclusion of the constant loops.   This isomorphism implies that the structure group of $\xi$ 
can be reduced to $U(n) \subset LU(n)$.  

The converse of this statement is immediate.  If the structure group of $\xi$  can be
reduced to $U(n) \subset LU(n)$ then $\xi \cong L\bc \otimes U(\xi) $.  If $\psi \subset \xi$ is the subbundle corresponding
to $L_+\bc \otimes U(\xi)$  then $\psi$ defines a Fourier decomposition of $\xi$. 
 \end{proof}

\med
Notice in the proof of this theorem we observed that if $\psi \subset \xi$ yields a Fourier decomposition  then
there is an isomorphism, $L\bc \otimes  (\psi \cap z\psi^\perp) \cong \xi$.  This says  that $\bc [z, z^{-1}] \otimes (\psi \cap
z\psi^\perp) $ can be viewed as a dense subbundle  of $\xi$  or equivalently,  there is a fiberwise dense   filtration,
$$
\cdots \subset z^{-k}\psi \subset z^{-(k+1)}\psi \subset \cdots \subset \xi.
$$
This is the fiberwise analogue of the Fourier decomposition of $L\bc^n$.

\med
\bf Remarks on the proofs of corollaries 2 - 4.   \rm As observed in the introduction, corollary 2 is simply a restatement of theorem 1. 
Note also that corollary 3 follows from our proof of theorem 1.   To prove corollary 4,  let
$\Omega\zeta $ be the  principal $\om U(n)$ bundle over the based loop space classified by
$\om f_\zeta : \om M \to \om BU(n) \simeq U(n)$.  The argument in the proof of theorem \ref{main} proves that if $L\zeta$ has a
Fourier decomposition determined by $\psi \subset L\zeta$,  there is an isomorphism of vector bundles whose structure group is $\om
U(n)$,
\begin{equation}\label{iso}
\rho : \om \bc \otimes (\psi \cap z\psi^\perp)_{|_{\om M}} \cong \om \zeta.
\end{equation}
 But, as shown in the proof of theorem \ref{main}, the  bundle $\psi \cap z\psi^{\perp}$ is classified by the composition
 
$$
\begin{CD}
 LM @>Lf_\zeta >> LBU(n) @>e>> BU(n)  
\end{CD}
$$
which  is homotopic to
the composition
$$
\begin{CD}
  LM @>e>> M @>f_\zeta >> BU(n).
\end{CD}
$$
This composition is null homotopic when restricted to the based loop space, $\om M$, so the restriction $(\psi \cap
z\psi^\perp)_{|_{\om M}}$ is a  trivial $n$-dimensional bundle.  Thus (\ref{iso}) yields a trivialization of $\om \zeta$ as vector
bundle with structure group $\om U(n)$.  This implies the classifying map $\om f : \om M \to U(n)$ is null homotopic, which is  the
statement of   corollary 4.

\bg
\section{Twisting a loop bundle by a gauge transformation}

\med
Let $\zeta \to M$ be a smooth $n$ - dimensional complex bundle classified by a map
$f_\zeta : M \to B\gln$.  Let $L\zeta \to LM$ be the loop bundle described above.  It is classified by  the loop  of
$f$,    $Lf : LM \to LB\gln$.   Let
$e  : LX
\to X$ be the evaluation map  and consider the pull -back bundle,
$e^*(\zeta) \to LM$.  In this section we describe how to twist the loop bundle $L\zeta$ by an automorphism of
the bundle $e^*(\zeta)$  and we prove theorem \ref{gauge} as stated in the introduction.   

Let $\cge$ be the gauge group of smooth bundle automorphisms of $e^*(\zeta)$.  An element $g \in \cge$ can be
viewed as a diagram
$$
\begin{CD}
e^*(\zeta)   @>g >>  e^*(\zeta) \\
@VVV   @VVV \\
LM  @>>= >  LM.
\end{CD}
$$

We now go about defining the loop bundle $L^{\cg}\zeta \to \cge \times LM  $ discussed in theorem 
\ref{gauge}.  We do this first by defining a classifying map $$L^\cg f_\zeta :  \cge \times LM  \to LB\gln.$$

Recall that $\cge$ can be described as the group of sections of the adjoint bundle of the associated
principal bundle to $e^*(\zeta)$.  This principal bundle is the pull back of the universal bundle $E\gln \to B\gln$  via
the composition $f_\zeta \circ e : LM \to B\gln$.   Let $Ad (\zeta)$ denote the corresponding fiber bundle:
$$
Ad (\zeta) = ((f_\zeta \circ e)^* (E\gln) \times_{Ad} \gln \la LM
$$ where the notation $\times_{Ad}$ refers to taking the orbit space of the
diagonal $\gln$ action, where $\gln$ acts freely on $(f_\zeta \circ e)^* (E\gln) $  and by conjugation
on $\gln$.   The following is standard (see \cite{atiyahbott} for a good reference).

\begin{proposition}  There is a natural isomorphism of groups,
$$
\cge \cong \Gamma(Ad (\zeta)),
$$
where the right hand side is the group of all smooth sections of the adjoint bundle.  The multiplication in $\Gamma
(Ad(\zeta))$ is the fiberwise product of sections.
\end{proposition}

Now for any group $G$, let $Ad(EG) = EG\times_{Ad}G$ be the corresponding adjoint bundle.
The following is well known.

\med
\begin{proposition}  There is a natural homotopy equivalence $h: Ad(EG) \to LBG$ making the following map of
fibrations homotopy commute:
$$\begin{CD}
G   @>h >\simeq >   \Omega BG \\
@VVV @VVV \\
Ad(EG)  @>h>\simeq >  LBG \\
@VVV   @VVe V \\
BG   @>> = >  BG
\end{CD}
$$ where $h : G \to \Omega BG $ is the usual homotopy equivalence, adjoint to the map $\Sigma G \to BG$ that
classifies the $G$ - bundle over $\Sigma G$ with the single clutching function given by the identity $G \to G$.
\end{proposition}

\med
By these two propositions we have the following composition:
\begin{multline}
\bar\Phi_\zeta : \cge \cong \Gamma (Ad(e^*\zeta)) \hk Map (LM,
Ad(e^*\zeta)) \\
 \xrightarrow{(f_\zeta \circ e)^*}   Map (LM, Ad(EG) )
\xrightarrow[\cong]{h} (L M, LBG).
\end{multline}
Notice that the map $\bar \Phi_\zeta : \cge \to Map(LM, LBG)$ is well defined up to homotopy.  Therefore the adjoint
map
$$
\Phi_\zeta : \cge \times LM \to LBG
$$
is well defined up to homotopy and hence determines an isomorphism class of loop bundle 
\begin{equation}\label{Phi}
L^\cg \zeta \to \cge
\times LM.
\end{equation}

We can describe the bundle $L^\cg \zeta$ in a fiberwise manner as follows.   Let $\gamma \in LM$.  
For $t \in S^1$, let $t\gamma$ be rotation of the loop $\gamma$ by $t$.  That is, $t\gamma (s) = \gamma (s+t)$.
   Now let $\tau \in \cge$.  Then $\tau (\gamma) : \zeta_{\gamma (0)} \to \zeta_{\gamma (0)}$ is an
isomorphism.   So for $t \in S^1$,  $\tau (t\gamma) : \zeta_{\gamma (t)} = \zeta_{t\gamma (0)} \to \zeta_{t\gamma (0)} =
\zeta_{\gamma (t)}$ is an isomorphism.  Finally, let $\tilde \gamma : \br \to M$ be the composition $\begin{CD} \tilde \gamma :
\br \to
\br/\bz @>\gamma >> M \end{CD}$.

 The fiber $L^\cg \zeta_{|_{(\tau, \gamma)}}$ is given by
\begin{equation}\label{fiber}
L^\cg \zeta_{|_{(\tau, \gamma)}} = \{\sigma \in \Gamma (\tilde \gamma^*(\zeta)) : \ \text{for each} \ t \in \br,
\ \sigma (t + 1) = \tau(t\gamma)( \sigma (t))\}.
\end{equation}

For $\tau \in \cge$, let $L^\tau \zeta \to LM$ be the loop bundle given by the restriction of 
$L^\cg\zeta$ to  $\{\tau\} \times LM$.  As can be seen by either considering the classifying map $\Phi_\zeta$ or
from the fiberwise description
\ref{fiber}, one notices that for the identity
$ id \in \cge$, $ L^{id}(\zeta) = L\zeta \to LM$.  Therefore we view the bundles $L^\tau \zeta$ as twistings
or deformations of the  loop bundle $L\zeta$ by the automorphism $\tau$.

\med

\bf Interpretation.  \rm  Each $\tau \in \cge$ is an automorphism of $e^*\zeta$  and therefore determines
a fiberwise action of the integers $\bz$ on $e^*\zeta \to LM$.  That is, for each loop $\gamma$ there is a
representation  of $\bz$  on $\zeta_{\gamma (0)}$ given by $\tau (\gamma)$.  In particular this induces a fiberwise
action of $\bz$ on $\tilde \gamma^*\zeta \to \br$,   where the representation of the fiber over $t \in \br$, which
is  $ \zeta_{\gamma (t)}$, is given by $\tau (t\gamma)$.  Of course there is the free action of
$\bz$ on
$\br$, and the fiber at each loop, $L^\tau\zeta_\gamma$,   consists of  the $\bz$-equivariant  sections of  
$\tilde \gamma^*(\zeta)$:
$$
(L^\tau \zeta)_\gamma  = \Gamma_\bz (\tilde \gamma^*\zeta).
$$
  When $\tau = id$,   the resulting $\bz$-action on each
$\tilde \gamma*\zeta$ is trivial  and  so the fiber of the \emph{untwisted} loop bundle, $L\zeta$, at $\gamma \in
LM$  consists of the \emph{invariant} sections of $\tilde \gamma^*(\zeta)$, which are the sections of
$\gamma^*(\zeta) \to \br/\bz$. 

\med

An important example of a deformation (twisting)  of the loop bundle comes from a connection $\alpha$ on $\zeta$.
For  such a connection, let $\tau_\alpha$ be the parallel transport operator.  We can view $\tau_\alpha$ as an
element of the gauge group $\cge$, by defining for every $\gamma \in LM$, 
$\tau_\alpha (\gamma ) : \zeta_{\gamma (0)} \to \zeta_{\gamma(0)}$ to be the holonomy operator
around the loop $\gamma$ determined by the connection $\alpha$.  We refer to the corresponding loop bundle
$L^{\tau_\alpha}\zeta$ as the \emph{holonomy  loop bundle} associated to the connection $\alpha$.  Notice
that since the space of connections is affine, and therefore connected, the subspace of $\cge$ determined by
connections on $\zeta$ is connected.  Therefore the isomorphism class of the holonomy bundle
$L^{\tau_\alpha}\zeta$ is independent of the  choice of  connection $\alpha$.

We now discuss the equivariance of the bundle $L^\cg \zeta$.  Consider the usual rotation action
on the loop space, $S^1 \times LM \to LM$, and the trivial  $S^1$ action on $\cge$.  We then take the diagonal
action of
$S^1$ on the product, $\cge \times LM$.

\med
\begin{proposition}  The loop bundle $L^\cg \zeta$ has an $\br$-action, yielding an action map of loop
bundles,
$$
\begin{CD}
\br \times_{\bz} L^\cg \zeta  @>>> L^\cg \zeta \\
@VVV    @VVV \\
\br/\bz \times (\cge \times LM) @>>> \cge \times LM.
\end{CD}
$$
\end{proposition}

\begin{proof}  To define the action we consider the  fiberwise description of $L^\cg \zeta$ given  above
(\ref{fiber}).   Let
$\sigma \in L^\cg\zeta_{(\tau, \gamma)}$, where $(\tau,\gamma) \in \cge \times LM$.   Let $t\in \br$.  Define
$t\sigma \in L^\cg\zeta_{(\tau, t\gamma)}$ by $t\sigma (s) = \sigma (s+t)$ for $s\in \br$.  It is immediate that
this defines an action that satisfies the required properties. 
\end{proof}

Notice that since the circle action on $\cge$ is trivial, this induces an $\br$-action on each twisted loop bundle,
$L^\tau \zeta$ for  $\tau \in \cge$.  In particular for $\tau = id$, this action descends to the usual circle
action $\br/\bz \times L\zeta \to L\zeta$.  For $\zeta = TM$, then $LTM = TLM$  is the tangent bundle of the
loop space  and this action is the tangential action induced by rotation of loops.  The twisted tangent loop bundles,
$L^\tau TM$,  do not in general have $S^1$-actions, but only $\br$-actions. 

Notice that besides the $\br$-action, each twisted bundle, $L^\tau \zeta$, is a module bundle over the ring of
functions, $L\bc$.   The action is given by pointwise multiplication of sections.  That is, if $f \in L\bc$, and
$\sigma \in  \{\Gamma (\tilde \gamma^*(\zeta)) : \ \text{for each} \ t \in \br,
\ \sigma (t + 1) = \tau(t\gamma)( \sigma (t))\} = L^\tau \zeta_{|_\gamma}$,   then 
\begin{equation}
f\cdot \sigma (t) = f(t)\sigma (t)
\end{equation}
is a section of $\tilde \gamma^*(\zeta) $ living in $L^\tau \zeta_{|_\gamma}$.   

\med
The following  result will complete the proof of theorem 5, as stated in the introduction.

\med
\begin{proposition}  Let $\tau_\alpha \in \cge$ be the parallel transport operator of a connection $\alpha$ on $\zeta$.
Then there is an  $\br$-equivariant  isomorphism of  $L\bc$-module loop bundles,
$$
j : e^*\zeta \otimes L\bc \to L^{\tau_\alpha}\zeta.
$$
  
\end{proposition}

\begin{proof}  We begin by defining the embedding $j : e^*\zeta \hk L^{\tau_\alpha}\zeta$.  For $\gamma \in LM$
and  $t \in \br$, let $\tau_\alpha^{[0,t]}(\gamma) : \zeta_{\gamma (0)} \to \zeta_{\gamma (t)}$ denote the  parallel
transport operator along the curve $\gamma$ between $0$ and $t$.  We then define $j : e^*\zeta \hk
L^{\tau_\alpha}\zeta$ in a fiberwise way as follows.  Let $v \in \zeta_{\gamma (0)} = e^*\zeta_{\gamma}$. 
Define $j (v) \in \Gamma(\tilde \gamma^*(\zeta))$ by $$j(v)(t) = \tau_\alpha^{[0,t]} (\gamma) (v) \in
\zeta_{\gamma (t)} =
\tilde \gamma^*(\zeta)_t.$$  By construction, $j(v)(t +1) = \tau_\alpha (t\gamma) (j(v)(t))$ and so $j(v) \in
L^{\tau_\alpha}(\zeta)_\gamma$.  Clearly this defines an embedding of bundles over $LM$,
$j : e^*\zeta \hk L^{\tau_\alpha}\zeta$.  Moreover, this embedding is $\br$-equivariant  where $\br$ acts
on $e^*(\zeta)$ by $t\cdot v = \tau_\alpha^{[0,t]} (\gamma)(v)$ for $v \in (e^*\zeta)_\gamma = 
\zeta_{\gamma (0)}.$
Using the $L\bc$-module structure, we now extend this embedding to:

$$
j : e^*\zeta \otimes L\bc \to L^{\tau_\alpha}\zeta.
$$
This is an $\br$-equivariant map of $L\bc$-module loop bundles.  We now only need to verify that it is an
isomorphism.  For this we define an inverse map on a  fiber  over a loop $\gamma \in LM$.   The fiber over the left
hand side,  $ e^*\zeta \otimes L\bc$ at $\gamma$, is  $\zeta_{\gamma (0)} \otimes L\bc  = 
L\zeta_{\gamma (0)} $.   The fiber over the right hand side,
$(L^{\tau_\alpha}\zeta)_\gamma$,  is given by the space of sections  $\{\sigma \in \Gamma (\tilde \gamma^*(\zeta)) : \ \text{for
each} \ t \in \br,
\ \sigma (t + 1) = \tau_\alpha(t\gamma)( \sigma (t))\}$  (see (\ref{fiber})).  We now define an inverse map
\begin{equation}\label{inverse}
\Phi_\gamma : (L^{\tau_\alpha}\zeta)_\gamma \la  
L\zeta_{\gamma (0)}.
\end{equation}

For   $t \in \br$, let $-\tilde \gamma_t : [0, t] \to M$ be given by $-\tilde \gamma_t (s) = \tilde \gamma (t-s)$.   Then define
$$
\Phi_\gamma (\sigma) (t) =  \tau_\alpha^{[0,t]}(-\tilde \gamma_t) (\sigma (t)) \in \zeta_{\gamma (0)},
$$
where, as before,  $\tau_\alpha^{[0,t]}(-\tilde\gamma_t) : \zeta_{\gamma (t)} \to \zeta_{\gamma (0)}$
is parallel transport along the path $-\tilde \gamma_t$  via the connection $\alpha$.   To check it is well defined, we
need to see that $\Phi_\gamma (\sigma) (t+1) = \Phi_\gamma (\sigma) (t) $ for all $t \in \br$.  This is true because
 $\sigma (t + 1) = \tau_\alpha(t\gamma)( \sigma (t)) $.  One now immediately sees that $\Phi_\gamma$ is inverse to 
$j_\gamma :  \zeta_{\gamma (0)} \otimes \lc \to (L^{\tau_\alpha}\zeta)_\gamma$.
\end{proof}

\med
\bf Remark.  \rm  Notice that this theorem says that the holonomy loop bundle $L^{\tau_\alpha}\zeta$ has a Fourier decomposition for any
connection $\alpha$.  On the other hand, Corollary \ref{description} of theorem \ref{main} says that the \sl untwisted \rm
loop bundle $L\zeta = L^{id}\zeta$ has a Fourier decomposition if and only if $L\zeta \cong \zeta \otimes \lc \cong L^{\tau_\alpha}\zeta$.
Viewing  the holonomy loop bundle $ L^{\tau_\alpha}\zeta$ as a \sl deformation \rm of $L\zeta$ by means of the bundle $L^\cg \zeta$ over
$\cg \times LM$, we see that $L\zeta$ has a Fourier decomposition if the parallel transport operator $\tau_\alpha$ lies in component of the identity
in $\cg$.

\bg


\begin{thebibliography}{99}
\bibitem{atiyahbott} M. Atiyah and R. Bott, \emph{The Yang-Mills equations over Riemann surfaces}, Phil. Trans. R. Soc. Lond. A
\bf 308, \rm  (1982), 523-615.

 \bibitem{morava} J. Morava, \emph{The tangent bundle of an almost-complex free loopspace},   Proceedings of  Stanford
workshop on equivariant homotopy theory:  Homology, Homotopy, and Applications \bf 3 \rm (2001), 407-415 

\bibitem{segpres}
 A. Pressley and G. Segal
\textbf{ Loop Groups},
  Oxford Math. Monographs, Clarendon Press
(1986).
 
\end{thebibliography}
\end{document}